\documentclass[12pt,a4paper]{article}
\usepackage[utf8x]{inputenc}
\usepackage[english]{babel}
\usepackage{amsmath}
\usepackage{amsfonts}
\usepackage{amssymb}
\usepackage{cite}
\usepackage[top=2cm, bottom=2.5cm, left=2.5cm, right=2cm]{geometry}
\usepackage{wasysym} %для прямых интегралов (немецкий стиль), нужно использовать \varint
\let\oldint\varint
\renewcommand\int{\oldint\limits}
\begin{document}

\begin{center}
\textbf{AN INVERSE PROBLEM FOR STURM–LIOUVILLE OPERATORS ON THE HALF-LINE HAVING BESSEL-TYPE SINGULARITY IN AN INTERIOR POINT}\\[0.3cm]

{\textbf{A.~Fedoseev}\\[0.6cm]}
\end{center}

\begin{abstract}
We study the inverse problem of recovering Sturm--Liouville operators on the half-line with a Bessel-type singularity inside the interval from the given Weyl function. The corresponding uniqueness theorem is proved, a constructive procedure for the solution of the inverse problem is provided, also necessary and sufficient conditions for the solvability of the inverse problem are obtained.\\


\end{abstract}

\section*{1. Introduction}

Consider the differential equation
\begin{equation}
\ell y = -y''+\Big(\frac{\nu_0}{(x-a)^2}+q(x)\Big)y=\lambda y,\quad x>0,
\label{initeq}
\end{equation}
on the half-line with a Bessel-type singularity at an interior point $a>0$. Here $q(x)$ is a complex-valued function and $\nu_0$ is a complex number. Let $\nu_0=\nu^2-1/4$ and to be definite, we assume that  $Re\,\nu>0,$ $\nu\neq 1,2,\ldots$ (the other cases require minor modifications). We also assume that $q(x)|x-~a|^{\min(0,1-2Re\,\nu)} \in L(0,T)$ for some $T>a$ and $q(x)\in L(T,\infty)$. We denote the class of such functions $q(x)$ by $W$.

The paper deals with the boundary value problem  ${\cal L=L}(q)$ for differential equation \eqref{initeq} with the Dirichlet boundary condition $y(0)=0$ and an additional \textit{matching conditions} near the singular point $x=a$. We consider in some sense arbitrary matching conditions with a transition matrix $A=[a_{jk}]_{j,k=1,2},$ that connects solutions of \eqref{initeq} in a neighbourhood of the singular point (for details, see Section 2). 
In the special case without singularities ($\nu_0=0$) the above-mentioned matching conditions correspond to the conditions
\begin{equation}
\left[\begin{array}{l}  y \\ y'  \end{array} \right] (a+0)=
A\left[\begin{array}{l} y \\ y'  \end{array} \right] (a-0).
\label{matchcondwithoutnu}
\end{equation}

The goal of the present paper is to study the inverse problem of reconstructing ${\cal L}$ from the given Weyl function.
Inverse problems for differential equations with singularities inside an 
interval are important in mathematics and its applications. A wide class of differential
equations with turning points can be reduced to equations with singularities. For example, inverse
problems for such equations occur in electronic engineering in designing heterogeneous transmission
channels with given characteristics \cite{yurko-freiling-99}. Boundary value problems with a discontinuity at an interior
point also arise in geophysical models of the Earth’s crust \cite{lapwood81}. Inverse problems for equations
with singularities and turning points are used in investigations of the discontinuous solutions of
some integrable nonlinear equations of mathematical physics \cite{const98}. Note that in different problems
of natural sciences, various matching conditions for solutions with the corresponding transition
matrices $A$ are used.

Direct and inverse problems for the classical Sturm-Liouville operators without singularities were studied fairly completely (see \cite{rundell97,marchenko77,levitan84,yurko-NOVA,mclaughlin86} and references therein). 
In the case of $\nu_0=0$, self-adjoint operators with the conditions of the form \eqref{matchcondwithoutnu} were studied in \cite{hald84,shep94,yurko-00-IT}.
The case in which a singularity lies at an endpoint of the interval (i.e., $a = 0$) was considered in \cite{stash53,yurko-92-DU,yurko-93-IP,yurko-95-MS}. 
% The inverse problem of the recovering equation \eqref{initeq} from the given spectral data was studied in \cite{yurko-02-DU}. 
In \cite{yurko-02-DU}, the equation \eqref{initeq} under the Robin boundary condition was considered and the inverse problem of recovering equation \eqref{initeq} from the given spectral data  under spectrum simplicity restriction was studied. 
The uniqueness theorem for the higher-order equation on the half-line having nonintegrable singularity was proved in \cite{tamkang11}.

Spectral analysis of the singular non-self-adjoint boundary value problem ${\cal L}$ on the half-line
with arbitrary matching conditions produces essential qualitative modifications in the analytic
techniques, especially in inverse spectral problems. 
For the boundary value problem with discontinuities, the behaviour of the spectrum is more
complicated than for classical Sturm–Liouville operators. In particular, the discrete spectrum
can be unbounded and partially lie on the continuous spectrum. 
% In contrast to the previous work \cite{yurko-02-DU} in this paper we will not require any conditions on the spectrum.
In contrast to the work \cite{yurko-02-DU} in this paper we investigate the inverse problem of recovering equation \eqref{initeq} under the Dirichlet boundary condition from the given Weyl function and we will not require any conditions on the spectrum.
The method of transformation operators applied in \cite{marchenko77,levitan84,hald84,shep94} to classical Sturm–Liouville operators is inconvenient for problem ${\cal L}$.
The analysis of the inverse problem in the present paper is based on the ideas of the method of the spectral mappings \cite{yurko-VSP}. In Section 2, we study spectral properties of the boundary value problem ${\cal L}$ and introduce for the problem ${\cal L}$ so-called Weyl function which is a generalization of the notion of the Weyl function for the classical Sturm-Liouville operators \cite{marchenko77,levitan84, yurko-NOVA}. In Sections 3 and 4, we investigate the inverse problem of recovering ${\cal L}$ from the given Weyl function. The uniqueness theorem is proved and a constructive procedure for the solution of the inverse problem is obtained in Section 3. Necessary and sufficient conditions for the solvability of the inverse problem are obtained in Section 4.

\section*{2. Spectral properties}

\medskip
Let $\lambda=\rho^2$ and $\mathrm{Im}\,\rho\geq0$. Consider the functions
\[
C_j(x,\lambda)=(x-a)^{\mu_j} \sum^\infty_{k=0} c_{jk}(\rho (x-a))^{2k} ,\; j=1,2,
\]
where $c_{10}c_{20}=(2\nu)^{-1}$,
\[
\mu_j=(-1)^j\nu+1/2,\, c_{jk}=(-1)^k c_{j0}
\,\Big( \prod^k_{s=1} ((2s+\mu_j)(2s+\mu_j-1)-\nu_0)\Big)^{-1}.
\]
Here and below $z^{\mu}=\exp(\mu(\mbox{ln}|z|+i\arg z))$, $\arg z \in (-\pi, \pi]$. If $x>a$ or $x<a$ then the functions
 $C_j(x,\lambda)$ are solutions of \eqref{initeq} with $q(x)\equiv 0$. Let the functions $s_j(x,\lambda)$, $j=1,2$ be solutions of the following integral equations for $x>a$ and $x<a$: 
\[
s_j(x,\lambda)=C_j(x,\lambda)+\int^x_{a} g(x,t,\lambda)q(t)s_j(t,\lambda)\,dt,
\]
where $g(x,t,\lambda)=C_1(t,\lambda)C_2(x,\lambda)-C_1(x,\lambda)C_2(t,\lambda)$. For each fixed $x$ the functions $s_j(x,\lambda)$ are entire in $\lambda$ of order 1/2 and form a fundamental system of solutions of \eqref{initeq}. Moreover 
\begin{equation}
\det[s_j^{(m-1)}(x,\lambda)]_{j,m=\overline{1,2}}\equiv 1. 
\label{detsj}        
\end{equation} 

Let $A=[a_{jk}]_{j,k=1,2}$, $\det A\ne 0$ be a given matrix with complex numbers $a_{jk}$. We introduce the functions  $\{\sigma_j(x,\lambda)\}_{j=1,2}$, $x\in J_{-}\cup J_{+}$, $J_{\pm}=\{\pm (x-a)>0\}$ by the formula
\begin{equation*}
\sigma_j(x,\lambda)=\left\{
\begin{array}{ll}
s_j(x,\lambda), & x\in J_{-},\\
\sum\limits_{k=1}^{2}a_{kj}s_k(x,\lambda),&  x\in J_{+}.
\end{array}\right.
\label{sigmajinsk}
\end{equation*}
The fundamental system of solutions $\{\sigma_j(x,\lambda)\}$ is used to match solutions in a neighborhood of the singular point  $x=a$. It follows from \eqref{detsj} that
\begin{equation}
\det[\sigma_j^{(m-1)}(x,\lambda)]_{j,m=1,2}\equiv
\left\{ \begin{array}{ll}
1,\quad & x\in J_{-}, \\ \det A,\quad & x\in J_{+}.
\end{array}\right. 
\label{detsigmaj}                
\end{equation}

We introduce numbers $\xi_{jk}$, $j,k=1,2$ by
\begin{equation}
\left[ \begin{array}{ll}
\xi_{11} & \xi_{12} \\ \xi_{21} & \xi_{22}
\end{array}\right]= \frac{1}{2\sin\pi\nu}
\left[ \begin{array}{ll}
-a_{11}e^{2\pi i\nu}+a_{22}e^{-2\pi i\nu} & -i(a_{11}e^{\pi i\nu}-a_{22}e^{-\pi i\nu}) \\
-i(a_{11}e^{\pi i\nu}-a_{22}e^{-\pi i\nu}) & a_{11}-a_{22} 
\end{array}\right] .    
\label{ksimatrix}
\end{equation}       
% The behaviour of the spectrum of the boundary value problem ${\cal L}$ depends on the numbers $\xi_{jk}$. To be
% definite, we consider the most important special case in which $|\xi_{jj}|>|\xi_{12}|>0$ and $a_{12}=0$. In contrast to the case of classical Sturm–Liouville operators, here the discrete spectrum is unbounded and essential qualitative modifications
% in the investigation of the direct and inverse problems are arised.
The behaviour of the spectrum of the boundary value problem ${\cal L}$ depends on the numbers $\xi_{jk}$. We consider the most difficult special case in which $|\xi_{jj}|>|\xi_{12}|>0$, since in this case the discrete spectrum is unbounded and essential qualitative modifications in the investigation of the direct and inverse problems are arised. 
To be definite we set $a_{12}=0$ (the other case requires minor modifications).

We set
\begin{equation}
\begin{array}{c}
\varphi_1(x,\lambda)=\sigma'_2(0,\lambda)\sigma_1(x,\lambda)-
\sigma'_1(0,\lambda)\sigma_2(x,\lambda),\\[3mm]
\varphi_2(x,\lambda)=\sigma_1(0,\lambda)\sigma_2(x,\lambda)-
\sigma_2(0,\lambda)\sigma_1(x,\lambda).
\end{array}
\label{phiinsigmasys}                                      
\end{equation}
The functions $\varphi_j(x,\lambda)$, $j=1,2,$ are solutions of equation \eqref{initeq} for $x\in J_{\pm}$, and satisfy the initial conditions 
\begin{equation*}
\varphi_j^{(m-1)}(0,\lambda)=\delta_{jm},\quad j,m=1,2,                
\label{phicond}
\end{equation*}
where $\delta_{jm}$ is the Kronecker delta. By virtue of \eqref{detsigmaj} and \eqref{phiinsigmasys},
\begin{equation}
\det[\varphi_j^{(m-1)}(x,\lambda)]_{j,m=1,2}\equiv
\left\{ \begin{array}{ll}
1,\quad & x\in J_{-}, \\ \det A,\quad & x\in J_{+}.
\end{array}\right.                                         
\label{detphij}
\end{equation}
Denote $[1]_a=1+O\Big(|\rho(x-a)|^{-1}+|\rho|^{-1}\Big)$, when $|\rho(x-a)|\geq1$, $|\rho|\to\infty$. The following Lemmas were proved in \cite{yurko-02-DU}.

\medskip
{\bf Lemma 1. }{\it If $|\rho(x-a)|\ge 1$, $m=0,1$, $|\rho|\to\infty$, then the asymptotic formulas
\begin{equation}
\varphi_j^{(m)}(x,\lambda)=\frac{1}{2}\Big((-i\rho)^{m-j+1}\exp(-i\rho x)[1]_a
+(i\rho)^{m-j+1}\exp(i\rho x)[1]_a\Big),\; x\in J_{-},          
\label{phijestx<a}
\end{equation}
\[
\varphi_j^{(m)}(x,\lambda)=\frac{1}{2}\Big((-i\rho)^{m-j+1}
(\xi_{12}\exp(-i\rho x)[1]_a+(-1)^{j-1}\xi_{22}\exp(-i\rho(x-2a))[1]_a)+
\]
\begin{equation}
+(i\rho)^{m-j+1}(\xi_{21}\exp(i\rho x)[1]_a+
(-1)^{j-1}\xi_{11}\exp(i\rho(x-2a))[1]_a)\Big),\; x\in J_{+},   
\label{phijestx>a}
\end{equation}
are valid, where the numbers $\xi_{jk}$ are defined in \eqref{ksimatrix}.}

By $\Pi_{+}$ we denote the $\lambda$-plane with the two-sided cut $\Pi_0$ along the ray $\Lambda_{+}:=\{\lambda:\, \lambda\ge 0\}$ and we set $\Pi:=\overline\Pi_{+}\setminus\{0\}$. Under the mapping $\rho\to\rho^2=\lambda$ the sets $\Pi_{+},\,\Pi_{0}$ and $\Pi$ correspond to the sets $\Omega_{+}=\{\rho:\;\mbox{Im}\,\rho >0\}$, $\Omega_0=\{\rho:\; \mbox{Im}\,\rho=0\}$ and $\Omega=\{\rho:\;\mbox{Im}\,\rho\ge 0,\,\rho\ne 0\}$ respectively.

We introduce the so-called \textit{discontinuous Jost solution} $e(x,\rho)$, $x\ge 0, \mbox{Im}\,\rho\ge 0$ of equation \eqref{initeq}, which can be expressed as
\begin{equation*}
e(x,\rho)=
\sum_{k=1}^{2} A_k(\rho)\sigma_k(x,\lambda),\quad x\in J_{\pm}, 
\label{jostinsigmak}
\end{equation*}
and satisfies the condition
\begin{equation}
\lim_{x\to\infty} e(x,\rho)\exp(-i\rho x)=1.             
\label{limjost}
\end{equation}

\medskip
{\bf Lemma 2. }{\it The function $e(x,\rho)$ has the following properties:

($i_1$) for each $x$ the functions $e^{(m)}(x,\rho)$, $m=0,1,$ are analytic in $\rho\in\Omega_{+}$  and continuous for $\rho\in\Omega$;

($i_2$) the following asymptotic formulas are valid for $|\rho(x-a)|\ge 1$, $m=0,1$, $|\rho|\to\infty$ as $|\rho|\to\infty$:
\begin{equation}
e^{(m)}(x,\rho)=(i\rho)^m \exp(i\rho x)[1]_a,
\quad x\in J_{+},\quad  \mbox{Im}\,\rho\ge 0,                 
\label{jostestx>a}
\end{equation}
\[
e^{(m)}(x,\rho)=(i\rho)^m (\det A)^{-1}\Big(\xi_{12}\exp(i\rho x)[1]_a+
\]
\begin{equation}
+(-1)^{m+1}\xi_{jj}\exp(i\rho(2a-x))[1]_a\Big),
\quad x\in J_{-},\; \rho\in\overline{S_{2-j}},\; j=1,2;      
\label{jostestx<a}
\end{equation}

($i_3$) for real values of $\rho\ne 0$, $x\in J_{\pm}$ the functions $e(x,\rho)$ and $e(x,-\rho)$ form a fundamental solution system of \eqref{initeq}; moreover,
\[
\langle e(x,\rho), e(x,-\rho)\rangle =
\left\{ \begin{array}{ll}
-2i\rho (\det A)^{-1}, & x\in J_{-}, \\
-2i\rho, &  x\in J_{+},
\end{array}\right.
\]
where $\langle y,z\rangle :=yz'-y'z.$}

Denote $S_{k_0}=\left\{\rho:\;\arg\rho\in\Big(\frac{k_0\pi}2,\frac{(k_0+1)\pi}2\Big)\right\}$, $k_0=0,1$ and 
\begin{equation}
\Delta(\rho)=e(0,\rho),\quad \mbox{Im}\,\rho\ge 0.    
\label{Delta}
\end{equation}
The function $\Delta(\rho)$ is called the \textit{characteristic function} of the boundary value problem ${\cal L}$.
It follows from Lemma 2 that the function $\Delta(\rho)$ is analytic for $\rho\in\Omega_{+}$ and continuous for $\rho\in\Omega$. Moreover,
\begin{equation}
\Delta(\rho)=(\det A)^{-1}
\Big([\xi_{12}]-[\xi_{jj}]\exp(2i\rho a)\Big), \quad
|\rho|\to\infty,\;\rho\in\overline{S_{2-j}},\;j=1,2,        
\label{Deltaest}
\end{equation}
where $[1]=1+O(\rho^{-1})$, $|\rho|\to\infty$.
For sufficiently large $|\rho|$ the function $\Delta(\rho)$ has countable set of zeros of the form
\begin{equation*}
\rho_k=\rho_k^{\pm}+O(k^{-1}), \quad k\to\pm\infty,           
\label{rhokestinrhok-+}
\end{equation*}
where the numbers $\rho_k^{\pm}=\frac{\pi}{a}(k+\theta_{\pm})$ are the zeros of the functions
\begin{equation*}
\Delta^{\pm}(\rho)=\xi_{12}-\xi_{jj}\exp(2i\rho a),\; \rho\in S_{2-j},\; j=1,2         
\label{Delta-+}
\end{equation*}
and
\begin{equation*}
\theta_{\pm}=-\frac{i}{2\pi}\ln\Big|\frac{\xi_{12}}{\xi_{jj}}\Big|
+\frac{1}{2\pi}\arg\Big(\frac{\xi_{12}}{\xi_{jj}}\Big)       
\label{theta-+}
\end{equation*}
($``-"$ corresponds to $j=1$ and $``+"$ corresponds to $j=2$). Obviously, $\mbox{Im}\,\theta_{\pm}>0$.
To be definite, we assume that $\arg\Big(\frac{\xi_{12}}{\xi_{jj}}\Big)\in [0,2\pi)$. We set
$\Lambda=\{\lambda=\rho^2: \rho\in\Omega,\;\Delta(\rho)=0\}$,
$\Lambda'=\{\lambda=\rho^2: \;\rho\in\Omega_{+},\,\Delta(\rho)=0\}$,
$\Lambda''=\{\lambda=\rho^2:\;\rho\in\Omega_0,\,\rho\ne 0,\,\Delta(\rho)
=0\}$. 
Obviously, $\Lambda=\Lambda'\cup\Lambda''$, $\Lambda'$ is a countable unbounded set and $\Lambda''$ is a bounded set.
Denote
\begin{equation}
\Phi(x,\lambda)=e(x,\rho)/ \Delta(\rho),\quad M(\lambda):=\Phi'(0,\lambda). 
\label{PhiandM}
\end{equation}
The function $\Phi(x,\lambda)$ is a solution of \eqref{initeq} and satisfies the conditions $\Phi(0,\lambda)=1$, $\Phi(x,\lambda)=O(\exp(i\rho x))$, $x\to\infty$, $\rho\in\Omega$. Function $\Phi(x,\lambda)$ is called the \textit{Weyl solution} for ${\cal L}$.
The function $M(\lambda)$ is called the {\it Weyl function} for ${\cal L}$. Let fixed matrix $A$ and number $\nu_0$ be given.

\medskip
\textbf{Problem 1.} \textit{Recover $q(x)$ by the given Weyl function $M(\lambda)$.}

\medskip
Here and in the sequel, $\varphi(x,\lambda):=\varphi_2(x,\lambda)$.
Clearly,
\begin{equation*}
M(\lambda)=e'(0,\rho)/ \Delta(\rho),                          
\label{Minjost}
\end{equation*}
\begin{equation}
\Phi(x,\lambda)=\varphi_1(x,\lambda)+M(\lambda)\varphi(x,\lambda).
\label{Phiinphi}
\end{equation}
Then by \eqref{detphij}, we have
\begin{equation}
\langle\Phi(x,\lambda),\varphi(x,\lambda)\rangle=\left\{
\begin{array}{ll}                                    
1, &\mathrm{for}\; x\in J_-,\\
\det A, &\mathrm{for}\; x\in J_+.
\end{array}\right.
\label{<phiPhi>}
\end{equation}
Relations \eqref{Delta} and \eqref{PhiandM}, in combination with Lemma 2, lead to the following assertion.

\medskip
{\bf Theorem 1. }{\it The Weyl function  $M(\lambda)$ is analytic in $\Pi_{+}\setminus\Lambda'$ and continuous
in $\Pi\setminus \Lambda$. The set of singularities of $M(\lambda)$ (as an analytic function) coincides with the set $\Lambda_0:=\Lambda_{+}\cup\Lambda$.}

\medskip
{\bf Definition 1. } The set of singularities of the Weyl function $M(\lambda)$ is called the \textit{spectrum} of ${\cal L}$ and is denoted by $sp\,{\cal L}$. The values of $\lambda$, for which equation \eqref{initeq} has a nontrivial solution of the form $y(x)=\sum_{k=1}^{2}a_k \sigma_k(x,\lambda)$, satisfying the conditions $y(0)=0,\; y(\infty)=0$ (i.e.,
 $\lim\limits_{x\to\infty} y(x)=0$), are called the \textit{eigenvalues} of ${\cal L}$, and the corresponding solutions are called \textit{eigenfunctions}.

\medskip
{\bf Theorem 2. } {\it If $\lambda>0$ then the boundary value problem ${\cal L}$ has no eigenvalues. Moreover, if $\lambda_0=\rho_0^2>0$ and $\Delta(\rho_0)=0$, then $\Delta(-\rho_0)\ne 0$.}

\medskip
{\bf Proof. } Suppose that  $\lambda_0 =\rho_0^2 >0$ is an eigenvalue and $y_0(x)$ is the corresponding eigenfunction. Since the functions $\{e(x,\rho_0), e(x,-\rho_0)\}$ form a fundamental solution system of \eqref{initeq} in $J_{\pm}$, we have $y_0(x)=C_1e(x,\rho_0)+C_2e(x,-\rho_0)$. As $x\to\infty$ we have $y_0(x)\sim 0$, $e(x,\pm\rho_0)\sim\exp(\pm i\rho_0 x)$. But this is possible only if $C_1=C_2=0$. Further, if $\lambda_0 =\rho_0^2 >0$ and $\Delta(\rho_0)=0$, then $0\ne\langle e(x,\rho_0)$, $e(x,-\rho_0)\rangle_{|x=0} =-e'(0,\rho_0)\Delta(-\rho_0)$, and consequently $\Delta(-\rho_0)\ne 0$. The proof is complete.
$\hfill\Box$

\medskip
{\bf Theorem 3. }{\it The set $\Lambda'=\{\lambda_k\}$ coincides with the set of nonzero eigenvalues of problem ${\cal L}$, and}
\begin{equation}
e(x,\rho_k)=\beta_k \varphi(x,\lambda_k),
\quad \beta_k \ne 0,\quad \lambda_k=\rho_k^2.         
\label{jostinphi}
\end{equation}

\medskip
{\bf Proof. } Let $\lambda_k\in\Lambda'$. 
It follows from \eqref{Delta} that $e(0,\rho_k)=\Delta(\rho_k)=0$, and relation \eqref{limjost}, implies that $\lim\limits_{x\to\infty} e(x,\rho_k)=0$. Therefore, $e(x,\rho_k)$ is an eigenfunction, and $\lambda_k=\rho_k^2$ is an eigenvalue. In addition, it follows from \eqref{PhiandM} and \eqref{<phiPhi>} that $\langle e(x,\rho_k),\varphi(x,\lambda_k)\rangle=0$, i.e., relation \eqref{jostinphi} holds.

Conversely, let $\lambda_k=\rho_k^2,\;\mbox{Im}\,\rho_k>0$ be an eigenvalue and $y_k(x)$ be an eigenvalue. Obviously, $y_k'(0)\ne 0$. Without loss of generality, we can assume that $y_k'(0)=1$. Whence it follows that $y_k(x)\equiv\varphi(x,\lambda_k).$ Since $\lim\limits_{x\to\infty} y_k(x)=0$, we have $y_k(x)=\beta_k^0 e(x,\rho_k)$, $\beta_k^0\ne 0$. This implies \eqref{jostinphi}. Therefore, $\Delta(\rho_k)=e(0,\rho_k)=0$, and $\varphi(x,\lambda_k)$ and
$e(x,\rho_k)$ are eigenfunctions. The proof is complete.
$\hfill\Box$

We set
\[
G_{\delta}:=
\{\rho:\;\mbox{Im}\,\rho\ge 0,\;|\rho-\rho_k|\ge\delta,\;\rho_k\in\Lambda\}.
\]
Using \eqref{Deltaest}, \eqref{PhiandM} and Lemma 2, we obtain for $|\lambda|\to\infty$, $\rho\in G_\delta\cap\overline S_{2-j}$, $j=1,2$ 
\begin{equation}
M(\lambda)=i\rho\Big(M_0^\pm(\lambda)+O\big(\frac1\rho\big)\Big),
\label{Mineq}
\end{equation}
\[
M_0^\pm(\lambda)=\frac{\xi_{12}+\xi_{jj}\exp(2i\rho a)}{\xi_{12}-\xi_{jj}\exp(2i\rho a)},
\]
where $``-"$ corresponds to $j=1$ and $``+"$ corresponds to $j=2$.

\section*{3. Solution of the inverse problem}

In this section we prove the uniqueness theorem for the solution of this inverse problem. For this purpose we use ideas of the contour integral method.
In the analysis of the inverse problem, along with ${\cal L}$ we consider a boundary value problem  $\tilde{\cal L}$ of the same form but with different coefficients  $\tilde q$. If a certain symbol $\gamma$ denotes an object related to ${\cal L}$, then the corresponding symbol $\tilde\gamma$ with tilde will denote the analogous object related to $\tilde{\cal L}$, and  $\widehat\gamma:=\gamma-\tilde\gamma$.

\medskip
{\bf Theorem 4. }{\it If $M(\lambda)=\widetilde M(\lambda)$, then $q(x)=\tilde q(x)$ a.e. for $x>0$. Thus, the Weyl function uniquely determines the boundary value problem $\cal L$. }

\medskip
{\bf Proof. } Define matrix $P(x,\lambda)=[P_{kj}(x,\lambda)]_{k,j=1,2}$ by formulas
\begin{equation}
\begin{gathered}
P_{k1}(x,\lambda)=-\frac1{\eta(x)}\Big(\varphi^{(k-1)}(x,\lambda)\tilde\Phi'(x,\lambda)-\Phi^{(k-1)}(x,\lambda)\tilde\varphi'(x,\lambda)\Big),\\
P_{k2}(x,\lambda)=-\frac1{\eta(x)}\Big(\Phi^{(k-1)}(x,\lambda)\tilde\varphi(x,\lambda)-\varphi^{(k-1)}(x,\lambda)\tilde\Phi(x,\lambda)\Big),
\end{gathered}
\label{Pk12}
\end{equation}
where $\eta(x)=1$ if $x\in J_{-}$, and $\eta(x)=\det A$ if $x\in J_{+}$.
In virtue of \eqref{<phiPhi>} it follows, that
\begin{equation}
\begin{gathered}
\varphi(x,\lambda)=P_{11}(x,\lambda)\tilde\varphi(x,\lambda)+P_{12}(x,\lambda)\tilde\varphi'(x,\lambda),\\
\Phi(x,\lambda)=P_{11}(x,\lambda)\tilde\Phi(x,\lambda)+P_{12}(x,\lambda)\tilde\Phi'(x,\lambda).
\end{gathered}
\label{phiandPhiinPk12}
\end{equation}
It follows from \eqref{PhiandM}, Lemma 1 and Lemma 2 that
\begin{equation}
\begin{gathered}
|\varphi^{(m)}(x,\lambda)|\leq C|\rho|^{m-1}|\exp (-i\rho x)|,\\
|\Phi^{(m)}(x,\lambda)|\leq C_\delta|\rho|^{m}|\exp (i\rho x)|,\quad \rho\in G_\delta,\\
|\varphi^{(m)}(x,\lambda)-\tilde\varphi^{(m)}(x,\lambda)|\leq C|\rho|^{m-2}|\exp (-i\rho x)|,\\
|\Phi^{(m)}(x,\lambda)-\tilde\Phi^{(m)}(x,\lambda)|\leq C_\delta|\rho|^{m-1}|\exp (i\rho x)|,\quad \rho\in G_\delta \cap \tilde G_\delta.
\end{gathered}
\label{phiandPhiest}
\end{equation}
for $x>0$ and $m=0,1$.
Using \eqref{<phiPhi>}, \eqref{Pk12} and \eqref{phiandPhiest} we obtain
\begin{equation}
P_{jk}(x,\lambda)-\delta_{jk}=O(\rho^{-1}),\quad j\le k; \qquad P_{21}(x,\lambda)=O(1),
\label{Pkest}
\end{equation}
for $x\ge 0$ and $\rho\in G_\delta$, $|\rho|\to\infty$.
By the assumptions we have $M(\lambda)=\widetilde M(\lambda)$. Then, from \eqref{Pk12} and \eqref{Phiinphi}, we conclude, that for each fixed $x$ functions $P_{jk}(x,\lambda)$ are entire functions of $\lambda$. Taking into account \eqref{Pkest}, we obtain $P_{11}(x,\lambda)\equiv 1$, $P_{12}(x,\lambda) \equiv 0$. Substituting these relations into \eqref{phiandPhiinPk12}, we find that $\varphi(x,\lambda)\equiv\tilde\varphi(x,\lambda)$, $\Phi(x,\lambda)\equiv\tilde\Phi(x,\lambda)$ for all $x$ and $\lambda$, whence it follows that ${\cal L}=\tilde{\cal L}$, which completes the proof of Theorem 4. 
$\hfill\Box$

\medskip
Now we are going to construct the solution of the inverse problem. We say that ${\cal L} \in V$, if $q(x)\in W$. The inverse problem we will be solved in the class $V$.

Let $\tilde {\cal L}={\cal L}(\tilde q)$ is chosen such, that $\widehat M(\lambda)=O(1)$ (it can always be done due to \eqref{Mineq}).
Denote
\begin{equation*}
D(x,\lambda,\mu)=\frac{1}{\eta(x)}\frac{\langle \varphi(x,\lambda),
\varphi(x,\mu)\rangle}{\lambda-\mu},\quad
r(x,\lambda,\mu)=D(x,\lambda,\mu)\widehat M(\mu).
\label{Ddef}
\end{equation*}
It follows from Lemma 1 that if $\lambda=\rho^2$, $\mu=\theta^2$, $0\le \mbox{Im}\,\rho\le C$, $0\le \mbox{Im}\,\theta\le C$, then for fixed $x\in J_{\pm}$ the following estimates are valid:
\begin{equation*}
|D(x,\lambda,\mu)|\le\frac{C}{|\rho||\theta|(|\rho\mp\theta|+1)},\quad
|\varphi(x,\lambda)|\le C,
\quad \pm \mbox{Re}\,\rho\,\mbox{Re}\,\theta\ge 0.
\label{Dest}
\end{equation*}
Functions $\tilde r$ and $\tilde D$ are defined similarly, but with $\tilde\varphi$ instead of $\varphi$.
Let us take $h>0$ such, that $\mbox{Im}\,\rho_k<h,\; \mbox{Im}\,\tilde\rho_k<h$ for each $\rho_k\in
\Lambda,\; \tilde\rho_k\in\tilde\Lambda.$ Let $\gamma=\{\lambda=u+iv:\;
u=(2h)^{-2}v^2-h^2\}$ is image of the set $\mbox{Im}\,\rho=h$ by the mapping
$\lambda=\rho^2$. We set $J_\gamma=\{\lambda:\;\lambda\not\in\gamma\,\cup\,\mbox{int}\,\gamma\}$. 

\medskip
{\bf Theorem 5. }{\it The following relations hold}
\begin{gather}
\label{maineq}
\tilde \varphi(x,\lambda)=\varphi(x,\lambda)-\frac1{2 \pi i}\int\limits_{\gamma} \tilde r(x,\lambda,\mu)\varphi(x,\mu)\,d\mu,\\
\label{mainrel}
r(x,\lambda,\mu)-\tilde r(x,\lambda,\mu)-\frac1{2 \pi i}\int\limits_{\gamma} \tilde r(x,\lambda,\xi)r(x,\xi,\mu)\,d\xi =0,\\
\label{mainPhieq}
\tilde \Phi(x,\lambda)=\Phi(x,\lambda)-\frac1{2 \pi i}\frac1{\eta(x)}\int\limits_{\gamma} \frac{\langle\tilde\Phi(x,\lambda),\tilde\varphi(x,\mu)\rangle}{\lambda-\mu}\widehat M(\mu)\varphi(x,\mu)\,d\mu,\quad \lambda\in J_\gamma.
\end{gather}
\textit{Here (and further when it is necessary) the integral is treated in the sense of principal value:} $\int_\gamma=\lim\limits_{R\to\infty}\int_{{\gamma}_R}$.

\medskip
{\bf Proof. }
We choose positive numbers $r_N=((N+\chi)\pi/a)^2$ so, that circles 
$\theta_N:=\{\lambda:\; |\lambda|=r_N\}$ are lying in $G_{\delta}$ for sufficient small $\delta>0$. We set $\theta_{N,0}=\{\lambda:\;|\lambda|\le r_N\}$, $\gamma_N=(\gamma\cap\theta_{N,0})\cup\{\lambda:\;|\lambda|=r_N, \;\lambda\in \mbox{int}\,\gamma\}$ (with counterclockwise orientation).
According to Cauchy's integral formula we have
\[
P_{1k}(x,\lambda)=\delta_{1k}+\frac1{2\pi i}\int\limits_{\gamma_N} \frac{P_{1k}(x,\mu)}{\lambda-\mu}\,d\mu+
\frac1{2\pi i}\int\limits_{\theta_N} \frac{P_{1k}(x,\mu)-\delta_{1k}}{\lambda-\mu}\,d\mu, \; \lambda\not\in\mathrm{int}\,\gamma_N.
\]
Using \eqref{Pkest}, we get
\[
\lim_{N\to\infty} \frac1{2\pi i}\int\limits_{\theta_N} \frac{P_{1k}(x,\mu)-\delta_{1k}}{\lambda-\mu}\,d\mu=0,
\]
and, therefore
\begin{equation}
P_{1k}(x,\lambda)=\delta_{1k}+\frac1{2\pi i}\int\limits_{\gamma} \frac{P_{1k}(x,\mu)}{\lambda-\mu}\,d\mu,\quad\lambda\in J_\gamma.
\label{P1keq}
\end{equation}

By virtue of \eqref{phiandPhiinPk12} and \eqref{P1keq}
\[
\varphi(x,\lambda)=\tilde \varphi(x,\lambda)+\frac1{2 \pi i}\int\limits_{\gamma} 
\frac{\tilde\varphi(x,\lambda)P_{11}(x,\mu)+\tilde\varphi'(x,\lambda)P_{12}(x,\mu)}{\lambda-\mu}\,d\mu,\;\lambda\in J_\gamma.
\]
Using \eqref{Pk12}, we obtain
\begin{gather*}
\varphi(x,\lambda)=\tilde \varphi(x,\lambda)-\frac1{2 \pi i}\frac1{\eta(x)}\int\limits_{\gamma} 
\bigg(\tilde\varphi(x,\lambda)\Big(\varphi(x,\mu)\tilde\Phi'(x,\mu)-\Phi(x,\mu)\tilde\varphi'(x,\mu)\Big)+\\
+\tilde\varphi'(x,\lambda)\Big(\Phi(x,\mu)\tilde\varphi(x,\mu)-\varphi(x,\mu)\tilde\Phi(x,\mu)\Big)\bigg)\,\frac{d\mu}{\lambda-\mu}.
\end{gather*}
Hence, and from \eqref{Phiinphi} it follows \eqref{maineq}, since terms with $\varphi(x,\mu)$ are equal to nil by the Cauchy's theorem.

Relations \eqref{mainrel} and \eqref{mainPhieq} can be derived in a similar way. Theorem 5 is proved.
$\hfill\Box$

\medskip
For each given $x\geq0$, relation \eqref{maineq} can be treated as a linear equation with respect to $\varphi(x,\lambda)$. We will call \eqref{maineq} the \textit{main equation} of inverse problem.

Consider the Banach space $C(\gamma)$ of continuous bounded functions $z(\lambda)$, $\lambda\in\gamma$ with the norm $\|z\|=\sup\limits_{\lambda\in\gamma}|z(\lambda)|$.

\medskip
{\bf Theorem 6. }{\it For each fixed $x\geq0$ the main equation \eqref{maineq} has the unique solution $\varphi(x,\lambda)\in C(\gamma)$.}

\medskip
{\bf Proof. }
For fixed $x\geq0$ we consider the following linear bounded operators in $C(\gamma)$:
\[
\tilde Az(\lambda)=z(\lambda)-\frac1{2 \pi i}\int\limits_{\gamma} \tilde r(x,\lambda,\mu)z(\mu)\,d\mu,\quad Az(\lambda)=z(\lambda)+\frac1{2 \pi i}\int\limits_{\gamma} r(x,\lambda,\mu)z(\mu)\,d\mu.
\]
Then
\[
\tilde AAz(\lambda)=z(\lambda)+\frac1{2 \pi i}\int\limits_{\gamma} 
\Big(r(x,\lambda,\mu)-\tilde r(x,\lambda,\mu)-
\frac1{2 \pi i}\int\limits_{\gamma} \tilde r(x,\lambda,\xi)r(x,\xi,\mu)\,d\xi\Big)z(\mu)\,d\mu.
\]
In view of \eqref{mainrel} this yields $\tilde AAz(\lambda)=z(\lambda)$, $z(\lambda)\in C(\gamma)$. Swapping $\cal L$ and $\tilde{ \cal  L}$, we similarly get $A\tilde Az(\lambda)=z(\lambda)$. Thus $\tilde AA=A\tilde A=E$, where $E$ is identity operator. Therefore, the operator $\tilde A$ has a bounded inverse operator and the main equation \eqref{maineq} can be uniquely solved for each $x\geq0$.
$ \hfill\Box$

Thus, we have obtained the following algorithm for the solution of the inverse problem.

\medskip
{\bf Algorithm 1. }{\it Let the function $M(\lambda)$ be given.}\\
1) {\it We choose $\tilde{\cal  L}\in V$.}\\
2) {\it We find $\varphi(x,\lambda)$ from the main equation \eqref{maineq}.}\\
3) {\it We construct $q(x)$ by the formula }
\[
q(x)=\lambda+\frac{\varphi''(x,\lambda)}{\varphi(x,\lambda)}-\frac{\nu_0}{(x-a)^2}.
\]

\section*{4. Necessary and sufficient conditions}
In this section we formulate necessary and sufficient conditions for the solvability of the inverse problem.
To simplify calculations we suppose that the model boundary value problem $\tilde {\cal L}={\cal L}(\tilde q)$ is chosen such, that for $|\rho|\to\infty$,
\begin{equation}
\widehat M(\lambda)=O\Big(\frac1{\rho^2}\Big).
\label{Mcondness4}
\end{equation}
In particular, if potentials have additional smootheness than condition \eqref{Mcondness4} is valid for any model boundary value problem $\tilde {\cal L}$.

We set
\begin{equation}
\varepsilon_0(x)=\frac1{2\pi i}\frac1{\eta(x)}\int\limits_\gamma \tilde\varphi(x,\mu)\varphi(x,\mu)\widehat M(\mu)\,d\mu,\quad
\varepsilon(x)=-2\varepsilon'_0(x),
\label{epsdef}
\end{equation}
\begin{equation*}
p(x):=\frac{\nu_0}{(x-a)^2}+q(x).
\label{p(x)}
\end{equation*}

\medskip
{\bf Theorem 7. }{\it The following relation holds}
\begin{equation}
\label{qinq}
q(x)=\tilde q(x)-\varepsilon(x).
\end{equation}

\medskip
{\bf Proof. } Differentiating \eqref{maineq} twice by $x$, and using \eqref{epsdef} and the relation
\[
\frac d{dx}\frac{\langle\varphi(x,\lambda),\varphi(x,\mu)\rangle}{\lambda-\mu}=\varphi(x,\lambda)\varphi(x,\mu),
\]
we get
\begin{equation}
\tilde \varphi'(x,\lambda)+\varepsilon_0(x)\tilde \varphi(x,\lambda)=
\varphi'(x,\lambda)-\frac1{2 \pi i}\int\limits_{\gamma} \tilde r(x,\lambda,\mu)\varphi'(x,\mu)\,d\mu,
\label{phi'eq}
\end{equation}
\begin{equation}
% \begin{array}{c}
\begin{gathered}
\tilde \varphi''(x,\lambda)=\varphi''(x,\lambda)-\frac1{2 \pi i}\int\limits_{\gamma} 
\tilde r(x,\lambda,\mu)\varphi''(x,\mu)\,d\mu-\\
-\frac1{2 \pi i}\frac1{\eta(x)}\int\limits_{\gamma} 
2\tilde\varphi(x,\lambda)\tilde\varphi(x,\mu)\widehat M(\mu)\varphi'(x,\mu)\,d\mu-\\
-\frac1{2 \pi i}\frac1{\eta(x)}\int\limits_{\gamma} 
\Big(\tilde\varphi(x,\lambda)\tilde\varphi(x,\mu)\Big)'\widehat M(\mu)\varphi(x,\mu)\,d\mu.
\end{gathered}
\label{phi''eq}
\end{equation}
Substituting in \eqref{phi''eq} derivative of the second order from the equation \eqref{initeq}, and then substituting $\varphi(x,\lambda)$, using \eqref{maineq} we get
\begin{equation*}
\begin{gathered}
\tilde p(x)\tilde\varphi(x,\lambda)=p(x)\tilde\varphi(x,\lambda)+
\frac1{2 \pi i}\frac1{\eta(x)}\int\limits_{\gamma} 
\langle\tilde\varphi(x,\lambda)\tilde\varphi(x,\mu)\rangle\widehat M(\mu)\varphi(x,\mu)\,d\mu+\\
+\frac1{2 \pi i}\frac1{\eta(x)}\int\limits_{\gamma} 
2\tilde\varphi(x,\lambda)\tilde\varphi(x,\mu)\widehat M(\mu)\varphi'(x,\mu)\,d\mu+\\
+\frac1{2 \pi i}\frac1{\eta(x)}\int\limits_{\gamma} 
\Big(\tilde\varphi(x,\lambda)\tilde\varphi(x,\mu)\Big)'\widehat M(\mu)\varphi(x,\mu)\,d\mu.
\end{gathered}
\end{equation*}
Hence we get \eqref{qinq}.
$\hfill\Box$

Let us formulate necessary and sufficient conditions for the solvability of the inverse problem.
Denote as $\mathbf{W}$ the set of functions $M(\lambda)$ such that \\
(i) $M(\lambda)$ is analytic in $\Pi_+$ excluding a countable set of poles $\Lambda'$ and continuous in $\Pi\setminus\Lambda$ ($\Lambda$ and $\Lambda'$ are depend of each function $M(\lambda)$).\\
(ii) holds \eqref{Mineq} when $|\lambda|\to\infty$.

\medskip
{\bf Theorem 8. }{\it The function $M(\lambda)\in \mathbf{W}$ is the Weyl function for ${\cal L} \in V$ then and only then, when the following conditions are hold:}\\
1) {\it (asymptotics) there exists $\tilde{\cal L}\in V$ such, that \eqref{Mcondness4} holds;}\\
2) {\it (S-condition) for each fixed $x\geq0$ equation \eqref{maineq} has a unique solution $\varphi(x,\lambda)\in C(\gamma)$;}\\
3) {\it $\varepsilon(x)\in W$, where the function $\varepsilon(x)$ is defined by \eqref{epsdef}.}\\
{\it Under these conditions the function $q(x)$ can be constructed by formula \eqref{qinq}.}

\medskip
The necessity of Theorem 8 is proved above. We prove the sufficiency. Let the function $M(\lambda)\in \mathbf{W}$ be given, and it satisfies conditions from Theorem 8, and let the function $\varphi(x,\lambda)$ be the solution of the main equation  \eqref{maineq}. Then \eqref{maineq} gives analytic continuation for $\varphi(x,\lambda)$ in the whole $\lambda$-plane, moreover  for each fixed $x\geq0$, the function $\varphi(x,\lambda)$ is entire in $\lambda$ of order 1/2. It can be shown that the functions $\varphi^{(\nu)}(x,\lambda)$, $\nu=0,1$, are absolutely continuous on the compacts if $|x-a|\geq\varepsilon$, for each fixed $\varepsilon>0$, and
\begin{equation}
|\varphi^{(\nu)}(x,\lambda)|\leq C|\rho|^{\nu-1} \exp(|\tau|x),\quad\lambda\in\gamma.
\label{phiestsuff}
\end{equation}
We construct the function $\Phi(x,\lambda)$ from the relations \eqref{mainPhieq}, and ${\cal L}={\cal L}(q)$ by formula \eqref{qinq}. It is clear, that ${\cal L}\in V$. 

\medskip
{\bf Lemma 3. }{\it The following relations hold}
\[
\ell\varphi(x,\lambda)=\lambda\varphi(x,\lambda),\quad\ell\Phi(x,\lambda)=\lambda\Phi(x,\lambda).
\]

\medskip
{\bf Proof. } Differentiating \eqref{maineq} twice by $x$ we get \eqref{phi'eq} and \eqref{phi''eq}. It follows from \eqref{phi''eq}, \eqref{maineq} and \eqref{qinq} it follows
\begin{equation}
\begin{gathered}
\tilde\ell\tilde\varphi(x,\lambda)=\ell\varphi(x,\lambda)-
\frac1{2 \pi i}\int\limits_{\gamma} \tilde r(x,\lambda,\mu)\ell\varphi(x,\mu)\,d\mu-\\
-\frac1{2 \pi i}\frac1{\eta(x)}\int\limits_{\gamma} 
\langle\tilde\varphi(x,\lambda),\tilde\varphi(x,\mu)\rangle\widehat M(\mu)\varphi(x,\mu)\,d\mu.
\end{gathered}
\label{lphieqsuff}
\end{equation}
Using \eqref{mainPhieq}, we similarly obtain
\begin{equation}
\begin{gathered}
\tilde\Phi'(x,\lambda)+\varepsilon_0(x)\tilde\Phi(x,\lambda)=\Phi'(x,\lambda)-\\
-\frac1{2 \pi i}\frac1{\eta(x)}\int\limits_{\gamma} 
\frac{\langle\tilde\Phi(x,\lambda),\tilde\varphi(x,\mu)\rangle}{\lambda-\mu}\widehat M(\mu)\varphi'(x,\mu)\,d\mu,
\end{gathered}
\label{Phieqsuff}
\end{equation}
\begin{equation}
\begin{gathered}
\tilde\ell\tilde\Phi(x,\lambda)=\ell\Phi(x,\lambda)-
\frac1{2 \pi i}\frac1{\eta(x)}\int\limits_{\gamma} 
\frac{\langle\tilde\Phi(x,\lambda),\tilde\varphi(x,\mu)\rangle}{\lambda-\mu}\widehat M(\mu)\ell\varphi(x,\mu)\,d\mu-\\
-\frac1{2 \pi i}\frac1{\eta(x)}\int\limits_{\gamma} 
\langle\tilde\Phi(x,\lambda),\tilde\varphi(x,\mu)\rangle\widehat M(\mu)\varphi(x,\mu)\,d\mu.
\end{gathered}
\label{lPhieqsuff}
\end{equation}
It follows from \eqref{lphieqsuff}, that
\begin{equation*}
\begin{gathered}
\lambda\tilde\varphi(x,\lambda)=\ell\varphi(x,\lambda)-
\frac1{2 \pi i}\int\limits_{\gamma} \tilde r(x,\lambda,\mu)\ell\varphi(x,\mu)\,d\mu-\\
-\frac1{2 \pi i}\int\limits_{\gamma} 
(\lambda-\mu)\tilde r(x,\lambda,\mu)\varphi(x,\mu)\,d\mu.
\end{gathered}
\end{equation*}
Taking \eqref{maineq} into account, we find, that for fixed $x\geq0$
\begin{equation}
\eta(x,\lambda)+\frac1{2 \pi i}\int\limits_{\gamma} \tilde r(x,\lambda,\mu)\eta(x,\mu)\,d\mu=0,\quad \lambda\in\gamma,
\label{etaeqsuff}
\end{equation}
where $\eta(x,\lambda)=\ell\varphi(x,\lambda)-\lambda\varphi(x,\lambda)$. According to \eqref{phiestsuff} we have for fixed $x\geq 0$ 
\begin{equation}
|\eta(x,\lambda)|\leq C|\rho|,\quad\lambda\in\gamma.
\label{etaest1}
\end{equation}
Using estimates \eqref{etaest1} and \eqref{etaeqsuff} we arrive to estimate $|\eta(x,\lambda)|\leq C$ for $\lambda\in\gamma$. By virtue of S-condition from Theorem 8 the homogeneous equation \eqref{etaeqsuff} have only trivial solution $\eta(x,\lambda)\equiv0$. Therefore,
\[
\ell\varphi(x,\lambda)=\lambda\varphi(x,\lambda).
\]
Hence and from \eqref{lPhieqsuff} together with \eqref{mainPhieq} we get $\ell\Phi(x,\lambda)=\lambda\Phi(x,\lambda)$.
$\hfill\Box$

\medskip
\textit{Let us continue the proof of Theorem 8}.
We set $x=0$ in \eqref{maineq}, \eqref{phi'eq}. Then we have
\begin{equation}
\varphi(0,\lambda)=\tilde\varphi(0,\lambda)=0,\quad \varphi'(0,\lambda)=\tilde\varphi(0,\lambda)=1.
\label{phicondsuff}
\end{equation}
Using \eqref{mainPhieq} and \eqref{Phieqsuff}, we calculate
\begin{equation}
\label{Phi0rel}
\Phi(0,\lambda)=\tilde\Phi(0,\lambda)=1,\quad
\Phi'(0,\lambda)=\tilde\Phi'(0,\lambda)+
\frac 1{2\pi i}\int\limits_\gamma\frac{\widehat M(\mu)}{\lambda-\mu}\,d\mu.
\end{equation}
We fix $\lambda\in J_\gamma$. 
From \eqref{mainPhieq}, \eqref{phicondsuff} taking into account the estimates  
\begin{equation*}
\begin{gathered}
|\tilde\varphi^{(m)}(x,\mu)|\le C|\theta|^{m-1}|\exp(-i\theta x)|,\quad \mu=\theta^2,\quad x\geq0,\quad m=0,1,\\[3mm]
|\tilde\Phi^{(m)}(x,\lambda)|\le C_\delta|\rho|^{m}|\exp(i\rho x)|,\quad x\geq0, \quad\rho\in G_\delta,
\end{gathered}
\end{equation*}
we get that
\[
\Phi(x,\lambda)=O(\exp(i\rho x+2h x)), x\to\infty.
\]
Hence, and since $\Phi(0,\lambda)=1$, it follows, that the function $\Phi(x,\lambda)$ is a Weyl solution.
Further, from \eqref{Phi0rel} we have
\[
\Phi'(0,\lambda)=\widetilde M(\lambda)+\frac1{2\pi i}\int\limits_\gamma\frac{\widehat M(\mu)}{\lambda-\mu}\,d\mu,
\]
According to Cauchy's integral formula 
\[
\widehat M(\lambda)=\frac1{2\pi i}\int\limits_{\gamma_N} \frac{\widehat M(\mu)}{\lambda-\mu}\,d\mu+
\frac1{2\pi i}\int\limits_{\theta_N} \frac{\widehat M(\mu)}{\lambda-\mu}\,d\mu, \; \lambda\not\in\mathrm{int}\,\gamma_N.
\]
Then for $N\to\infty$ we get
\[
\widehat M(\lambda)=\frac1{2\pi i}\int\limits_{\gamma} \frac{\widehat M(\mu)}{\lambda-\mu}\,d\mu,\quad\lambda\in J_\gamma.
\]
Therefore, $\Phi'(0,\lambda)=\widetilde M(\lambda)+\widehat M(\lambda)=M(\lambda)$, 
i.e. $M(\lambda)$ is the Weyl function for ${\cal L}$.
$\hfill\Box$

{\bf Remark 1.} Similar results for the Robin boundary condition $y'(0)-hy(0)=0$ are obtained in \cite{fedoseev-izvsgu-12}.

{\bf Acknowledgment.}  This research  was supported in part by Grants 
10-01-00099 and 10-01-92001-NSC of Russian Foundation for Basic Research
and Taiwan National Science Council.

\begin{minipage}[b]{0.49\linewidth}
\hspace{0.1cm}
\end{minipage}
\hfill
\begin{minipage}[b]{0.49\linewidth}
{
\small
Alexey Fedoseev

Department of Mathematics

Saratov State University

Astrakhanskaya 83, Saratov 410012, Russia

E-mail: fedoseev\_ae@mail.ru

}
\end{minipage}


\begin{thebibliography}{99}
\bibitem{yurko-freiling-99} 
	Freiling G and Yurko V A 1999 Reconstructing parameters of a medium from 
    incomplete spectral information {\it Results Math.} {\bf 35} 228-49
\bibitem{lapwood81} 	
	Lapwood F R and Usami T 1981 {\it Free Oscillations of the Earth}
    (Cambridge: Cambridge University Press)
\bibitem{const98} 		
	Constantin A 1998 On the inverse spectral problem for the Camassa-Holm equation
	{\it J. Funct. Anal.} {\bf 155} 352-63
\bibitem{rundell97}
	Chadan K, Colton D, P\"{a}iv\"{a}rinta L and Rundell W 1997 {\it An Introduction to Inverse Scattering 
	and Inverse Spectral Problems (SIAM Monographs on Mathematical Modeling and Computation)} (Philadelphia: SIAM)
\bibitem{marchenko77} 	
	Marchenko V A 1977 {\it Sturm-Liouville Operators and Their Applications} (Kiev: Naukova Dumka)\\ 
	Marchenko V A 1986 {\it Sturm-Liouville Operators and Their Applications} (Basle: Birkh\"auser)	(Engl. Transl.) 
\bibitem{levitan84}
	Levitan B M 1984 {\it Inverse Sturm-Liouville problems} (Moscow: Nauka) (Engl. Transl.)\\
	Levitan B M 1987 {\it Sturm-Liouville problems} (Utrecht: VNU Science)
\bibitem{yurko-NOVA} 	
	Freiling G and Yurko V A 2001 {\it Inverse Sturm-Liouville Problems
    and their Applications} (New York: NOVA Science)
\bibitem{mclaughlin86}    
    McLaughlin J R 1986 Analytical methods for recovering coefficients in differential equations from spectral data
    {\it SIAM Rev.} {\bf 28} 53-72
\bibitem{hald84}	
	Hald O H 1984 Discontinuous inverse eigenvalue problems
    {\it Commun. Pure Appl. Math.} {\bf 37} 539-77
\bibitem{shep94}
	Shepelsky D G 1994 The inverse problem of reconstruction of the medium`s conductivity in a class of discontinuous and
    increasing functions {\it Adv. Sov. Math.} {\bf 19} 209-231 
\bibitem{yurko-00-IT}    
    Yurko V A 2000 Integral transforms connected with  discontinuous boundary value problems
    {\it Integral Transforms and Special Functions} {\bf 10} 141-64
\bibitem{stash53}
	Stashevskaya V V 1953 On inverse problems of spectral analysis for a certain class of differential equations 
	{\it Dokl. Akad. Nauk SSSR} {\bf 93} 409-12	
\bibitem{yurko-92-DU}    
	Yurko V A 1992 Inverse problem for differential equations with a singularity {\it Differ. Uravneniya} {\bf 28} 1355-62\\
	Yurko V A 1992 {\it Differ. Equ.} {\bf 28} 1100-7 (Engl. Transl.)
\bibitem{yurko-93-IP}	
	Yurko V A 1993 On higher-order differential operators with a singular point
	{\it Inverse Problems} {\bf 9} 495-502
\bibitem{yurko-95-MS}	
	Yurko V A 1995 On higher-order differential operators with a singularity  {\it Matem. Sbornik} {\bf 186} 133-60\\
	Yurko V A 1995 {\it Sbornik: Mathematics} {\bf 186} 901-28 (Engl. Transl.)
\bibitem{yurko-02-DU}    
	Yurko V A 2002 On recovering singular non-selfadjoint differential operators with a singularity inside the interval {\it Differ. Uravneniya} {\bf 38} 645-59\\
	Yurko V A 2002 {\it Differ. Equ.} {\bf 38} 678-94 (Engl. Transl.)
\bibitem{tamkang11} 
	Fedoseev A 2011 Inverse problems for differential equations on the half-line having a singularity in an interior point			
	{\it Tamkang J. Math.} {\bf 42} 343-54
\bibitem{yurko-VSP} 	
	Yurko V A 2002 {\it Method of Spectral Mappings in the Inverse
    Problem Theory (Inverse and Ill-posed Problems Series)} (Utrecht: VSP)
\bibitem{fedoseev-izvsgu-12} 	
	Fedoseev A E 2012 Inverse problem for Sturm–-Liouville operator on the half-line having 
	nonintegrable singularity in an interior point {\it Izv. Saratov. Univ. Mat. Mekh. Inform.} {\bf 12:4} 49-55 (Russian)
    
\end{thebibliography}
\end{document}